\documentstyle{article}
\def\R{\ifmmode{\rm I\mkern-3.1mu
R\mkern1mu}\else{\rm I\kern-.18em  R\hskip1pt\ 
}\fi\relax}  
\font\bigreek=Symbol at 16pt
\font\smgreek=Symbol at 10pt
\def\equi{\longleftrightarrow}
\def\a{\alpha} 
\def\b{\beta} 
\def\g{\gamma}
\def\G{\Gamma}
\def\D{\triangle}
\def\t{\tau}
\def\d{\delta}  
\def\th{\theta} 
\def\l{\lambda}
\def\L{\Lambda}
\def\n{\nu}
\def\di{\diamond}

\def\ph{\phi}

\def\m{\mu}
\def\s{\sigma}
\def\sou{\overline}
\def\so{\underline} 
\def\O{\Omega}
\def\f{\rightarrow}
\def\q{\forall}
  
\def\v{\vdash}
 
\def\p{\succ}
 
\def\mats{\ifmmode{ {\hbox{\bigreek s}} }\else{ 
{\bigreek s} }\fi\relax}
\def\matsin{\ifmmode{ {\hbox{\smgreek s}} }\else{ 
{\smgreek s} }\fi\relax}
\def\matt{\ifmmode{ {\hbox{\bigreek t}} }\else{ 
{\bigreek t} }\fi\relax}
\def\mattin{\ifmmode{ {\hbox{\smgreek t}} }\else{ 
{\smgreek t} }\fi\relax}

\newtheorem{theo}{Theorem}[section]
\newtheorem{lemma}{Lemma}[section]

\newtheorem{corollary}{Corollary}[section]

\begin{document} 

\begin{center} 
\LARGE\bf 
Mixed Logic and Storage Operators\\[1cm]
\end{center}

\begin{center} 
\bf 
Karim NOUR\\
\rm
LAMA - Equipe de Logique, Universit\'e de Chamb\'ery\\
73376 Le Bourget du Lac\\
e-mail nour@univ-savoie.fr\\[1cm]
\end{center}

\begin{abstract} 

In 1990 J-L. Krivine introduced the notion of storage operators. They are
$\l$-terms which simulate call-by-value in the call-by-name strategy and they can be used in
order to modelize assignment instructions. J-L. Krivine has shown that there is a very
simple second order type in $AF2$ type system for storage operators using G\H{o}del
translation of classical to intuitionistic logic. \\
In order to modelize the control operators, J-L. Krivine has extended the system $AF2$
to the classical logic. In his system the property of the unicity of integers
representation is lost, but he has shown that storage operators typable in the
system $AF2$ can be used to find the values of classical integers. \\  
In this paper, we present a new classical type system based on a logical system called
mixed logic. We prove that in this system we can characterize, by types, the storage
operators and the control operators. We present also a similar result in the M.
Parigot's $\l \m$-calculus.

\end{abstract}

\section{Introduction}

In 1990, J.L. Krivine introduced the notion of storage operators (see [4]). They are closed
$\l$-terms which allow, for a given data type (the type of integers, for example), to simulate in
$\l$-calculus the "call by value" in a context of a "call by name" (the head reduction) and they can
be used in order to modelize assignment instructions. J.L. Krivine has shown that the formula $\q x
\{ N$*$[x] \f \neg\neg N[x] \}$ is a specification for storage operators for Church integers : where
$N[x]$ is the type of integers in $AF2$ type system, and the operation $*$ is the simple  G\H{o}del
translation from classical to intuitionistic logic which associates to every formula $F$ the
formula $F$* obtained by replacing in $F$ every atomic formula by its negation (see [3]).\\

The latter result suggests many questions :
\begin{itemize}
\item Why do we need a G\H{o}del translation ?
\item Why do we need the type $N$*$[x]$ which characterize a class larger than
integers ? \end{itemize}

In order to modelize the control operators, J-L. Krivine has extended the system $AF2$
to the classical logic (see [6]). His method is very simple : it consists of adding a
new constant, denoted by $C$, with the declaration $C : \q X \{ \neg \neg X \f X \}$
which axiomatizes classical logic over intuitionistic logic. For the constant $C$, he
adds a new reduction rule : $(C t t_1 ... t_n) \f (t \quad \l x (x \quad t_1 ... t_n))$ which is a
particular case of a rule given by Felleisen for control operator (see [1]). In this
system the property of the unicity of integers representation is lost, but J-L. Krivine
has shown that storage operators typable in the intuitionistic system $AF2$ can be used to find the
values of classical integers \footnote{The idea of using storage operators in classical logic is
due to M. Parigot (see [19])}(see [6]). \\

The latter result suggests also many questions :
\begin{itemize}
\item What is the relation between classical integers and the type $N$*$[x]$ ?
\item Why do we need intuitionistic logic to modelize the assignment
instruction and classical logic to modelize the control operators ? 
\end{itemize}

In this paper, we present a new classical type system based on a logical system called
mixed logic. This system allows essentially to distinguish between
classical proofs and intuitionistic proofs. We prove that, in this system, we can characterize, by
types, the storage operators and the control operators. This results give some answers
to the previous questions.\\

We present at the end (without proof) a similar result in the M. Parigot's $\l
\m$-calculus.\\

\bf Acknowledgement.\/ \rm We wish to thank J.L. Krivine, and C. Paulin for
helpful discussions. We don't forget the numerous corrections and suggestions from R. David and N.
Bernard.

\section{Pure and typed $\l$-calculus}

\begin{itemize}
\item Let $t,u,u_1,...,u_n$ be $\l$-terms, the application of $t$ to $u$ is denoted by $(t)u$. In
the same way we write $(t)u_1...u_n$ instead of $(...((t)u_1)...)u_n$.
\item $Fv(t)$ is the set of free variables of a $\l$-term $t$. 
\item The $\b$-reduction (resp. $\b$-equivalence) relation is denoted by $u \f\sb{\b} v$ (resp. $u
\simeq\sb{\b} v$).  
\item The notation $\s(t)$ represents the result of the simultaneous substitution $\s$
to the free variables of $t$ after a suitable renaming of the bounded variables of $t$.    
\item We denote by $(u)^n v$ the $\l$-term $(u)...(u)v$ where $u$ occurs $n$ times, and $\sou{u}$
the sequence of $\l$-terms $u_1,...,u_n$. If $\sou{u} = u_1,...,u_n$ $n \geq 0$, we denote by
$(t)\sou{u}$ the $\l$-term $(t)u_1...u_n$.   
\item Let us recall that a $\l$-term $t$ either has a head redex [i.e. $t=\l x_1 ...\l
x_n (\l x u) v v_1 ... v_m$, the head redex being $(\l x u) v$], or is in head normal form [i.e.
$t=\l x_1 ...\l x_n (x) v_1 ... v_m$]. The notation $u \p v$ means that $v$ is obtained from $u$
by some head reductions.  If $u \p v$, we denote by $h(u,v)$ the length of the head reduction
between $u$ and $v$. 
\end{itemize} 

\begin{lemma} (see[3]) \\
1) If $u \p v$, then, for any substitution $\s$, $\s(u) \p \s(v)$, and $h(\s(u),\s(v))$=h(u,v). \\ 
2) If $u \p v$, then, for every sequence of $\l$-terms $\sou{w}$, there is a $w$, such that
$(u)\sou{w} \p w$, $(v)\sou{w} \p w$, and $h((u)\sou{w},w)=h((v)\sou{w},w)+h(u,v)$.
\end{lemma}

\bf Remark. \rm Lemma 2.1 shows that to make the head reduction of $\s(u)$ (resp. of $(u)\sou{w}$)
it is equivalent - same result, and same number of steps - to make some steps in the head
reduction of $u$, and after make the head reduction of $\s(v)$ (resp. of $(v)\sou{w}$). $\Box$

\begin{itemize}
\item The types will be formulas of second order predicate logic over a given language. The logical
connectives are $\perp$ (for absurd), $\f$, and $\q$. There are individual (or first order)
variables denoted by $x,y,z,...,$ and predicate (or second order) variables denoted by
$X,Y,Z,....$  
\item We do not suppose that the language has a special constant for equality.
Instead, we define the formula $u=v$ (where $u,v$ are terms) to be $\q Y(Y(u) \f Y(v))$ where $Y$
is a unary predicate variable. Such a formula will be called an equation. We denote by $a \approx
b$, if $a=b$ is a consequence of a set of equations.
\item The formula $F_1 \f (F_2 \f(...\f (F_n \f G)...))$ is
also denoted by $F_1,F_2,...,F_n \f G$. For every formula $A$, we denote by $\neg A$ the
formula $A \f \perp$. If $\sou{v} = v_1,...,v_n$ is a sequence of variables, we denote by
$\q \sou{v} A$ the formula $\q v_1...\q v_n A$.  
\item Let $t$ be a $\l$-term, $A$ a type, $\G = x_1 : A_1 ,..., x_n : A_n$ a context, and $E$ a
set of equations. We define by means of the following rules the notion "$t$ is of
type $A$ in $\G$ with respect to $E$" ; this notion is denoted by  $\G\v_{AF2} t:A$ :
\end{itemize}     
\begin{itemize} 
\item[] (1) $\G\v_{AF2} x_i:A_i$ $1\leq i\leq n$. 
\item[] (2) If $\G,x:A \v_{AF2} t:B$, then $\G\v_{AF2} \l xt:A \f B$.  
\item[] (3) If $\G\v_{AF2} u:A \f B$, and $\G\v_{AF2} v:A$, then $\G\v_{AF2} (u)v:B$.
\item[] (4) If $\G\v_{AF2} t:A$, and $x$ is not free in $\G$, then $\G\v_{AF2} t:\q xA$. 
\item[] (5) If $\G\v_{AF2} t:\q xA$, then, for every term $u$, $\G\v_{AF2} t:A[u/x]$.   
\item[] (6) If $\G\v_{AF2} t:A$, and $X$ is not free in $\G$, then $\G\v_{AF2} t:\q XA$. 
\item[] (7) If $\G\v_{AF2} t:\q XA$, then, for every formulas $G$, $\G\v_{AF2} t:A[G/X]$. 
\item[] (8) If $\G\v_{AF2} t:A[u/x]$, and $u \approx v$, then $\G\v_{AF2}
t:A[v/x]$. \end{itemize} 

This typed $\l$-calculus system is called $AF2$ (for Arithm\'etique Fonctionnelle du
second ordre).
 
\begin{theo} (see [2]) The $AF2$ type system has the following properties :\\
1) Type is preserved during reduction.\\
2) Typable $\l$-terms are strongly normalizable.
\end{theo}

We present now a syntaxical property of system $AF2$ that we will use afterwards.

\begin{theo} (see [8])  If in the typing we go from $\G\v_{AF2} t:A$ to $\G\v_{AF2} t:B$,
then we may assume that we begin by the $\q$-elimination rules, then by the equationnal
rule, and finally  by the $\q$-introduction rules.
\end{theo}

\begin{itemize}
\item We define on the set of types the two binary relations $\lhd$ and $\approx$ as the least
reflexive and transitive binary relations such that :  
\begin{itemize} 
\item[] - $\q xA \lhd A[u/x]$, if $u$ is a term of language ;
\item[] - $\q XA \lhd A[F/X]$, if $F$ is a formula of language ; 
\item[] - $A \approx B$ if and only if $A=C[u/x]$, $B=C[v/x]$, and $u \approx v$. 
\end{itemize} 
\end{itemize} 

\section{Pure and typed $\l C$-calculus}

\subsection{The $C2$ type system}

We present in this section the J-L. Krivine's classical type system.

\begin{itemize}
\item We add a constant $C$ to the pure $\l$-calculus and we denote by $\L C$ the set of new
terms also called $\l C$-terms. We consider the following rules of reduction, called
rules of head $C$-reduction.
\begin{itemize} 
\item[] 1) $(\l x u) t t_1 ... t_n \f (u[t / x]) t_1 ... t_n$ for every $u, t, t_1,...,t_n \in
\L C$. 
\item[] 2) $(C) t t_1 ... t_n \f (t) \l x (x)t_1 ... t_n$ for every $ t, t_1,...,t_n
\in \L C$, $x$ being a $\l$-variable not appearing in $t_1,...,t_n$.
\end{itemize}
\end{itemize} 

\begin{itemize} 
\item For any $\l C$-terms $t,t'$, we shall write $t \p_C t'$ if $t'$ is obtained from $t$ by
applying these rules finitely many times. We say that $t'$ is obtained from $t$ by head
$C$-reduction.
\item A $\l C$-term $t$ is said $\b$-normal if and only if $t$ does not contain a
$\b$-redex.
\item A $\l C$-term $t$ is said $C$-solvable if and only if  $t \p_C (f)t_1,...,t_n$ where
$f$ is a variable.
\end{itemize}

It is easy to prove that : if $t \p_C t'$, then, for any substitution $\s$, $\s (t) \p_C \s
(t')$.

\begin{itemize}
\item We add to the $AF2$ type system the new following rule : 
\begin{center} 
(0) $\G \v C : \q X \{ \neg \neg X \f X \}$
\end{center}
This rule axiomatizes the classical logic over the intuitionistic logic.
We call $C2$ the new type system, and we write $\G \v_{C2} t : A$ if $t$ is of type $A$ in the
context $\G$.
\end{itemize}

It is clear that  $\G \v_{C2} t : A$ if and only if $\G , C : \q X \{ \neg \neg X \f X \} \v_{AF2} t :
A$.

\begin{theo} (see [6]) \\
1) If $\G \v_{C2} t:A$, and $t \f_{\b} t'$, then $\G \v_{C2} t':A$.\\
2) If $\G \v_{C2} t:\perp$, and $t \p_C t'$, then $\G \v_{C2} t':\perp$.\\
3) If $A$ is an atomic type, and $\G \v_{C2} t:A$, then $t$ is $C$-solvable.
\end{theo}

\subsection{The $M2$ type system}

In this section, we present the system $M2$. This system allows essentialy to
distinguish between classical proofs and intuitionistic proofs \\

We assume that for every integer $n$, there is a countable set of special $n$-ary second order
variables denoted by $X_C,Y_C,Z_C$...., and called classical variables.\\ 

Let $X$ be an $n$-ary predicate variable or predicate symbol. A type $A$ is said to be
ending with $X$ if and only if $A$ is obtained by the following rules :  
\begin{itemize}  
\item[] - $X(t_1,...,t_n)$ ends with $X$; 
\item[] - If $B$ ends with $X$, then $A \f B$ ends with $X$ for every type $A$ ;
\item[] - If $A$ ends with $X$, then $\q vA$ ends with $X$ for every variable $v$.\\
\end{itemize}
 
A type $A$ is said to be a classical type if and only if $A$ ends with $\perp$ or a
classical variable.\\

We add to the $AF2$ type system the new following rules :
\begin{itemize} 
\item[] (0$'$) $\G\v C : \q X_C \{ \neg \neg X_C \f X_C \}$
\item[] (6$'$) If $\G\v t:A$, and $X_C$ has no free occurence in $\G$, then $\G\v t: \q
X_C A$.  
\item[] (7$'$) If $\G\v t: \q X_C A$, and $G$ is a classical type, then
$\G\v t:A[G/ X_C]$.  
\end{itemize} 
 
We call $M2$ the new type system, and we write $\G \v_{M2} t:A$ if $t$ is of type $A$ in the
context $\G$.\\

We extend the definition of $\lhd$ by : $\q X_C A \lhd A[G / X_C]$ if $G$ is a classical type.

\begin{lemma} If $A$ is a classical type and $A \lhd B$ (or $A \approx B$), then $B$ is a
classical type.
\end{lemma}
\bf Proof \rm Easy. $\Box$

\subsection{The logical properties of $M2$}

We denote by $LAF2$, $LC2$, and $LM2$ the underlying logic systems of respectively $AF2$, $C2$,
and $M2$ type systems.\\

With each classical variable $X_C$, we associate a special variable $X^{\ast}$ of $AF2$
having the same arity as $X_C$. For each formula $A$ of $LM2$, we define the formula
$A$* of $LAF2$ in the following way :  \begin{itemize} 
\item[] - If $A=D(t_1,...,t_n)$ where $D$ is a predicate symbol or a predicate
variable, then $A$*=$A$ ;  
\item[] - If $A=X_C(t_1,...,t_n)$, then $A$*$=\neg X^{\ast}(t_1,...,t_n)$  ;  
\item[] - If $A=B \f C$, then $A$*$=B$*$ \f C$* ; 
\item[] - If $A=\q xB$, then $A$*=$\q xB$*.
\item[] - If $A=\q XB$, then $A$*=$\q XB$*.
\item[] - If $A=\q X_C B$, then $A$*=$\q X^{\ast} B$*.
\end{itemize} 

$A$* is called the G\H{o}del translation of $A$. 

\begin{lemma} If $G$ is a classical type of $LM2$, then $\v_{LAF2} \neg \neg G$*$ \equi G$*.
\end{lemma}
\bf Proof \rm It is easy to prove that $\v_{LAF2} G$*$ \f \neg \neg G$*.\\
We prove $\v_{LAF2} \neg \neg G$*$ \f G$* by induction on $G$.
\begin{itemize} 
\item[] - If $G = \perp$, then $G$*=$\perp$, and $\v_{LAF2} ((\perp \f \perp) \f \perp) \f
\perp$.
\item[] - If $G = X_C (t_1,...,t_n)$, then $G$*=$\neg X^{\ast}(t_1,...,t_n)$, and $\v_{LAF2} \neg
\neg \neg X^{\ast}(t_1,...,t_n) \f \neg X^{\ast}(t_1,...,t_n)$.
\item[] - If $G = A \f B$, then $B$ is a classical type and $G$* = $A$* $\f$ $B$*. By the induction
hypothesis, we have $\v_{LAF2} \neg \neg B$*$ \f B$*. Since $\v_{LAF2} \neg \neg
(A$*$\f B$*) $\f$ $(\neg \neg A$*$\f \neg \neg B$*), we check easily that $\v_{LAF2}
\neg \neg (A$* $\f B$*) $\f (A$* $\f B$*).   
\item[] - If $G = \q vG'$ where $v=x$ or $v=X$,
then $G'$ is a classical type and $G$*=$\q vG'$*. By the induction hypothesis, we have $\v_{LAF2}
\neg \neg G'$*$ \f G'$*. Since $\v_{LAF2} \neg \neg \q vG'$* $\f$ $\q v \neg \neg G'$*, we check
easily that $\v_{LAF2} \neg \neg \q vG'$* $\f \q vG'$*. 
\item[] - If $G = \q X_C G'$, then $G'$ is a classical type and $G$*=$\q X^{\ast} G'$*. By the
induction hypothesis, we have $\v_{LAF2} \neg \neg G'$*$ \f G'$*. Since $\v_{LAF2} \neg \neg \q
X^{\ast} G'$* $\f$ $\q X^{\ast} \neg \neg G'$*, we check easily that $\v_{LAF2} \neg
\neg \q X^{\ast} G'$* $\f \q X^{\ast} G'$*. $\Box$
\end{itemize}

\begin{lemma} Let $A,G$ be formulas of $LM2$, $t$ a term, $x$ a first order variable, and $X$ a
second order variable. We have :\\ 
1) $(A[t/x])$*$= A$*$[t/x]$.\\
2) $(A[G/X])$*$=A$*$[G$*$/X]$.
\end{lemma}
\bf Proof \rm By induction on $A$. $\Box$

\begin{lemma} Let $A$ be a formula of $LM2$, $G$ a classical type, and $X_C$ a
classical variable.\\ 
$\v_{LAF2} (A[G/X_C])$*$ \equi A$*$[\neg G$*$/X_C]$.
\end{lemma}
\bf Proof \rm By induction on $A$.
\begin{itemize} 
\item[] - If $A = D(t_1,...,t_n)$ where $D$ is a predicate variable or a predicate symbol, then
$A$*=$A$, and $\v_{LAF2} A \equi  A$. 
\item[] - If $A = X_C (t_1,...,t_n)$, then $A$*=$\neg X^{\ast}(t_1,...,t_n)$, and, by Lemma 3.2,
$\v_{LAF2} \neg \neg G$*$ \equi G$*. 
\item[] - If $A = B \f C$, then $A$* = $B$* $\f$ $C$*. By the induction hypothesis, we
have $\v_{LAF2} (B[G/X_C])$*$ \equi B$*$[\neg G$*$/X_C]$ and $\v_{LAF2}
(C[G/X_C])$*$ \equi C$*$[\neg G$*$/X_C]$. Therefore $\v_{LAF2} \{ (B[G/X_C])$*$
\f (B[G/X_C])$*$\} \equi \{ B$*$[\neg G$*$/X_C] \f C$*$[\neg G$*$/X_C] \}$.   
\item[] - If $A = \q vA'$, where $v=x$ or $v=X$, then $A$*=$\q vA'$*. By the induction hypothesis,
we have $\v_{LAF2} (A'[G/X_C])$*$ \equi A'$*$[\neg G$*$/X_C]$. Therefore $\v_{LAF2}
(\q vA'[G/X_C])$*$ \equi \q vA'$*$[\neg G$*$/X_C]$. 
\item[] - If $A = \q Y_C A'$, then $A$*=$\q Y^{\ast} A'$*. By the induction hypothesis, we
have $\v_{LAF2} (A'[G/X_C])$* 
$\equi A'$*$[\neg G$*$/X_C]$. Therefore $\v_{LAF2}
(\q Y_C A'[G/X_C])$*$ \equi$ $ (\q Y_C A')$*$[\neg G$*$/X_C]$. $\Box$
\end{itemize}

\begin{theo} If $A_1,...,A_n \v_{LM2} A$, then $A_1$*$,...,A_n$* $\v_{LAF2} A$*.
\end{theo}
\bf Proof \rm By induction on the proof of $A$ and using Lemmas 3.2, 3.3, and 3.4. $\Box$

\begin{corollary} Let $A,A_1,...,A_n$ be formulas of $LAF2$.\\ 
$A_1,...,A_n \v_{LM2} A$ if and only if  $A_1,...,A_n \v_{LAF2} A$.  
\end{corollary}
\bf Proof \rm We use Theorem 3.2. $\Box$\\

With each  predicate variable $X$ of $C2$, we associate a classical variable $X_C$ having the
same arity as $X$. For each formula $A$ of $LC2$, we define the formula $A^C$ of $M2$
in the following way :  \begin{itemize} 
\item[] - If $A=D(t_1,...,t_n)$ where $D$ is a constant symbol, then
$A^C=A$ ;  
\item[] - If $A=X(t_1,...,t_n)$ where $X$ is a predicate symbol, then
$A^C=X_C(t_1,...,t_n)$  ;   
\item[] - If $A=B \f C$, then $A^C=B^C \f C^C$ ; 
\item[] - If $A=\q xB$, then $A^C=\q xB^C$ ;
\item[] - If $A=\q XB$, then $A^C=\q X_CB^C$.
\end{itemize} 
$A^C$ is called the classical translation of $A$.

\begin{theo} Let $A_1,...,A_n,A$ be formulas of $LC2$.\\
 $A_1,...,A_n \v_{LC2} A$ if and only if $A_1^C,...,A_n^C \v_{LM2} A^C$.   
\end{theo}
\bf Proof \rm By induction on the proof of $A$. $\Box$

\section{Properties of $M2$ type system}

By corollary 3.1, we have that a formula is provable in system $LAF2$ if and only if it
is provable in system $LC2$. This resultat is not longer valid if we decorate
the demonstrations by terms. We will give some conditions on the formulas
in order to obtain such a result.\\

We define two sets of types of $AF2$ type system : $\O^+$ (set of $\q$-positive types), and $\O^-$
(set of $\q$-negative types) in the following way :
\begin{itemize} 
\item[] - If $A$ is an atomic type, then $A \in \O^+$, and $A \in \O^-$ ;
\item[] - If $T \in \O^+$, and $T' \in \O^-$, then, $T' \f T \in \O^+$, and $T \f T' \in \O^-$ ; 
\item[] - If $T \in \O^+$, then $\q x T \in \O^+$ ;
\item[] - If $T \in \O^-$, then $\q x T \in \O^-$ ;
\item[] - If $T \in \O^+$, then $\q X T \in \O^+$ ;
\item[] - If $T \in \O^-$, and $X$ has no free occurence in $T$, then $\q X T \in \O^-$.
\end{itemize} 

\begin{lemma}
1) If $A \in \O^+$ (resp. $A \in \O^-$) and $A \approx B$, then $B \in \O^+$ (resp.
$B \in \O^-$).\\
2) If $A \in \O^-$  and $A \lhd B \f C$, then $B \in \O^+$ and $C \in \O^-$.
\end{lemma} 
\bf Proof \rm Easy. $\Box$

\begin{theo} Let $A_1,...,A_n$ be $\q$-negative types, $A$ a $\q$-positive type of
$AF2$ which does not end with $\perp$, $B_1,...,B_m$ classical types, and $t$ a $\b$-normal
$\l C$-term.\\ If $\G = x_1:A_1,...,x_n:A_n,y_1:B_1,...,y_m:B_m \v_{M2} t:A$, then
$t$ is a normal $\l$-term, and $x_1:A_1,...,x_n:A_n \v_{AF2} t:A$.  
\end{theo} 
\bf Proof \rm We argue by induction on $t$.
\begin{itemize} 
\item[] - If $t$ is a variable, we have two cases :
\begin{itemize}
\item[] - If $t=x_i$ $1 \leq i \leq n$, this is clear.
\item[] - If $t=y_j$ $1 \leq j \leq m$, then $A=\q \sou{v} B$ where $B_j \lhd B'_j$ and $B'_j
\approx B$. Therefore, by Lemma 3.1, $A$ is a classical type. A contradiction.
\end{itemize}
\item[] - If $t=\l x u$, then $\G,x:E \v_{M2} u:F$, and $A=\q \sou{v}( E' \f F')$ where $E \approx
E'$,  $F \approx F'$ and $\sou{v}$ does not appear in $\G$. First, by Lemma 4.1,  $E \in \O^-$
and $F \in \O^+$, and then, by the induction hypothesis, $u$ is a normal $\l$-term, and
$x_1:A_1,...,x_n:A_n,x:E \v_{AF2} u:F$. Therefore $t$ is a normal $\l$-term, and
$x_1:A_1,...,x_n:A_n \v_{AF2} t:A$. 
\item[] - If $t=(x)u_1 ... u_r$ $r \geq 1$, we have two cases :
\begin{itemize}
\item[] - If $t=x_i$ $1 \leq i \leq n$, then $A_i \lhd B_1 \f C_1$, $C'_i \lhd B_{i+1} \f
C_{i+1}$ $1 \leq i \leq r-1$, $C'_r \lhd D$, $A = \q vD'$, where $C'_i \approx C_i$ $1
\leq i \leq r$, $D' \approx D$, and $\G \v_{M2} u_i:B_i$ $1 \leq i \leq r$. Since $A_i$ is a
$\q$-negative types, we prove (by induction and using Lemma 4.1) that for all $1 \leq i \leq r$
$B_i$ is a $\q$-positive types. By the induction  hypothesis we have $u_i$ is a normal $\l$-term,
and $x_1:A_1,...,x_n:A_n \v_{AF2} u_i:B_i$. Therefore $t$ is a normal $\l$-term, and
$x_1:A_1,...,x_n:A_n \v_{AF2} t:A$.
\item[] - If $t=y_j$ $1 \leq j \leq m$, then $B_j \lhd B_1 \f
C_1$, $C'_i \lhd B_{i+1} \f C_{i+1}$ $1 \leq i \leq r-1$, $C'_r \lhd D$, $A = \q vD'$,
where $C'_i \approx C_i$ $1 \leq i \leq r$, $D' \approx D$, and $\G \v_{M2} u_i:B_i$ $1 \leq i
\leq r$. Therefore, by Lemma 3.1, $A$ is a classical type. A contradiction.
\end{itemize} 
\item[] - If $t=(C)uu_1 ... u_r$ $r \geq 0$, then there is a classical type $E$ such that $\G \v_{M2}
u:\neg \neg E$, $E \lhd B_1 \f C_1$, $C'_i \lhd B_{i+1} \f C_{i+1}$ $1 \leq i \leq r-1$, $C'_r \lhd
D$, $A = \q vD'$, where $C'_i \approx C_i$ $1 \leq i \leq r$, $D' \approx D$, and $\G \v_{M2}
u_i:B_i$ $1 \leq i \leq r$. Therefore, by Lemma 3.1, $A$ is a classical type. A contradiction.
$\Box$
\end{itemize}

\begin{corollary} Let $A$ be a $\q$-positive type of $AF2$ and $t$ a $\b$-normal $\l C$-term. \\
If $\v_{M2} t:A$, then $t$ is a normal $\l$-term, and $\v_{AF2} t:A$. 
\end{corollary}
\bf Proof \rm We use Theorem 4.1. $\Box$\\

As for relation betwen the systems $C2$ and $M2$, we have the following result.

\begin{theo} Let $A_1,...,A_n,A$ be types of $C2$, and $t$ a $\l C$-term. \\
$A_1,...,A_n \v_{C2} t:A$ if and only if $A_1^C,...,A_n^C \v_{M2} t:A^C$. 
\end{theo}
\bf Proof \rm By induction on the typing of $t$. $\Box$

\section{The integers}

\begin{itemize} 
\item Each data type can be defined by a second order formula. For example, the type of integers
is the formula : $N[x]= \q X \{ X(0), \q y(X(y) \f X(sy)) \f X(x) \}$ where $X$ is a unary
predicate variable, $0$ is a constant symbol for zero, and $s$ is a unary function symbol for
successor. The formula $N[x]$ means semantically that $x$ is an integer if and only if $x$
belongs to each set $X$ containing $0$ and closed under the successor function $s$.\\ 
The $\l$-term $\so{0} = \l x \l fx$ is of type $N[0]$ and represents zero.\\ 
The $\l$-term $\so{s} = \l n\l x\l
f(f)((n)x)f$ is of type $\q y(N[y] \f N[s(y)])$ and represents the successor function.   
\item A set of equations $E$ is said to be adequate with the type of integers if and only if :
\begin{itemize} 
\item[] - $s(a) \not \approx 0$ ;
\item[] - If $s(a) \approx s(b)$ , then so is $a \approx b$.
\end{itemize}
In the rest of the paper, we assume that all sets of equations are adequate with the type of
integers.
\item For each integer $n$, we define the Church integer $\so{n}$ by $\so{n} = \l x\l f(f)^n x$.
\end{itemize}

\subsection{The integers in $AF2$}

The system $AF2$ has the property of the unicity of integers representation. 

\begin{theo} (see [2]) Let $n$ be an integer. If $\v_{AF2} t :N[s^n (0)]$, then $t \simeq\sb{\b}
\so{n}$. \end{theo} 

The propositional trace $N=\q X \{ X,(X \f X) \f X \}$ of $N[x]$ also defines the integers.

\begin{theo} (see [2])
If $\v_{AF2} t :N$, then, for a certain $n$, $t \simeq\sb{\b} \so{n}$.
\end{theo} 

\bf Remark\/ \rm A very important property of data type is the following (we express it for the type
of integers) : in order to get a program for a function $f : N \f N$ it is sufficient to prove $\v
\q x ( N[x] \f N[f(x)] )$. For example a proof of $\v \q x ( N[x] \f N[p(x)] )$ from the equations
$p(0)=0$, $p(s(x))=x$ gives a $\l$-term for the predecessor in Church intergers (see [2]). $\Box$
 
\subsection{The integers in $C2$}

The situation in system $C2$ is more complex. In fact, in this system
the property of unicity of integers representation is lost and we have only one
operational characterization of these integers. \\

Let $n$ be an integer. A classical integer of value $n$ is a closed $\l C$-term $\th_n$ such that
$\v_{C2} \th_n :N[s^n(0)]$.

\begin{theo} (see [6] and [12])
Let $n$ be an integer, and $\th_n$ a classical integer of value $n$.
\begin{itemize}
\item[] - if $n=0$, then, for every distinct variables $x,g,y$ : $(\th_n) x g y \p_C (x) y$ ;
\item[] - if $n \not = 0$, then there is $m \geq 1$ and a mapping  $I
: \{0,...,m \} \f N$, such that for every distinct variables $x,g,x_0,x_1,...,x_m$ :
\begin{itemize}
\item[] $(\th_n) x g x_0 \p_C (g) t_1 x_{r_0}$ ; 
\item[] $(t_i) x_i \p_C (g) t_{i+1} x_{r_i}$ $1\leq i\leq m$ ; 
\item[] $(t_m) x_m \p_C (x) x_{r_m}$ ;
\end{itemize}
where $I(0)=n$, $I(r_m)=0$, and $I(i+1)=I(r_i)-1$ $0\leq i\leq m-1$.
\end{itemize}
\end{theo}

We will generalize this result.\\

Let $O$ be a particular unary predicate symbol. The typed system $C2_O$ is the typed
system $C2$ where we replace the rules (2) and (7) by : 
\begin{itemize}
\item [] $(2_O)$ If $\G,x:A \v_{C2_O} t:B$, $A$ and $B$ are not ending with $O$, then $\G
\v_{C2_O} \l xt:A \f B$.  
\item [] $(7_O)$ If $\G\v_{C2_O} t:\q X A$, and $G$ is not ending with $O$, then $\G \v_{C2_O}
t:A[G/X]$. \end{itemize}

We define on the types of $C2_O$ a binary relation $\lhd_O$ as the least reflexive and transitive
binary relation such that :
\begin{itemize}
\item [] $\q xA \lhd_O A[u/x]$ if $u$ is a term of language ;
\item []  $\q XA \lhd_O A[G/X]$ if $G$ is a type which is not ending with $O$.
\end{itemize} 

\begin{lemma} a) If $\G \v_{C2_O} t:\perp$, and $t \p_C t'$, then $\G \v_{C2_O} t':\perp$.\\
b) If $\G \v_{C2_O} t:A$, and $A$ is an atomic type, then $t$ is $C$-solvable.
\end{lemma}
\bf Proof \rm
a)  It is enough to do the proof for one step of reduction. We have two cases :
\begin{itemize}
\item [] - If $t=(\l xu)vv_1...v_m$, then $t'=(u[v/x])v_1...v_m$, $\G,x:F \v_{C2_O} u:G$,  $F$ and
$G$ are not ending with $O$,  $G'\lhd_O F_1 \f G_1$, $G'_j \lhd_O F_{j+1} \f G_{j+1}$ $1 \leq j
\leq m-1$, $G_m \approx \perp$, $G_j \approx G'_j$ $1 \leq j \leq m-1$, $\G \v_{C2_O} v:F$, and
$\G \v_{C2_O} v_j:F_j$ $1 \leq j \leq m$. It is easy to check that $\G \v_{C2_O} u[v/x]:G$, then
$\G \v_{C2_O} t':\perp$. 
\item [] - If $t=(C)vv_1...v_m$, then $t'=(v)\l x(x)v_1...v_m$, and there is a type $A$ which is
not ending with $O$ such that : $A'\lhd_O F_1 \f G_1$, $G'_j \lhd_O F_{j+1} \f G_{j+1}$ $1 \leq j
\leq m-1$, $G_m \approx \perp$, $A \approx A'$, $G_j \approx G'_j$ $1 \leq j \leq m$, $\G \v_{C2_O}
v:\neg \neg A$, and $\G \v_{C2_O} v_j:F_j$ $1 \leq j \leq m$. It is easy to check that
$\G,x:A \v_{C2_O} (x)v_1...v_m:\perp$, but $A$ is not ending with $O$, then $\G \v_{C2_O}
\l x(x)v_1...v_m:\neg A$, and  $\G \v_{C2_O} t':\perp$. 
\end{itemize}   
b) Indeed, a typing of $C2_O$ may be seen as a typing of $C2$. $\Box$   

\begin{lemma}
a) If $\G \v_{C2_O} t:O(a)$, and $t \p_C t'$, then $t=t'$.\\
b) If $\G=y_1:A_1,...,y_n:A_n,x_1:O(a_1),...,x_m:O(a_m) \v_{C2_O} t:O(a)$, and all $A_i$ $1 \leq i\leq
n$ are not ending with $O$, then $t$ is one of $x_i$, and $a_i \approx a$  $1 \leq i \leq n$.
\end{lemma}
\bf Proof \rm
a) It is enough to do the proof for one step of reduction. We have two cases :
\begin{itemize}
\item [] - If $t=(\l xu)vv_1...v_m$, then $t'=(u[v/x])v_1...v_m$, $\G,x:F \v_{C2_O} u:G$, $F$ and
$G$ are not ending with $O$,  $G'\lhd_O F_1 \f G_1$, $G'_j \lhd_O F_{j+1} \f G_{j+1}$ $1 \leq j
\leq m-1$, $G_m \approx O(a)$, $G_j \approx G'_j$ $1 \leq j \leq m-1$, $\G \v_{C2_O} v:F$, and $\G
\v_{C2_O} v_j:F_j$ $1 \leq j \leq m$. Therefore $G_j$ $1 \leq j \leq m$ is not ending with $O$,
which is impossible since $G_m \approx O(a)$.  
\item [] - If $t=(C)vv_1...v_m$, then $t'=(v)\l x(x)v_1...v_m$, and there is a type $A$ which is
not ending with $O$ such that : $A'\lhd_O F_1 \f G_1$, $G'_j \lhd_O F_{j+1} \f G_{j+1}$ $1 \leq j
\leq m-1$, $G_m \approx O(a)$, $A \approx A'$, $G_j \approx G'_j$ $1 \leq j \leq m$, $\G \v_{C2_O}
v:\neg \neg A$, and $\G \v_{C2_O} v_j:F_j$ $1 \leq j \leq m$. $A$ is not ending with $O$, therefore
$G_j$ $1 \leq j \leq m$ is not ending with $O$, which is impossible since $G_m \approx O(a)$. 
\end{itemize}
b) By Lemma 5.1, we have $t \p_C (f)t_1...t_r$, and, by a), $t=(f)t_1...t_r$. Therefore $\G
\v_{C2_O} (f)t_1...t_r:O(a)$.  
\begin{itemize}
\item [] - If $f=x_i$ $1 \leq i \leq m$, then $r=0$, $t=x_i$, and $O(a_i) \approx O(a)$, then $a_i
\approx a$.
\item [] - If $f=y_j$ $1 \leq j \leq k$, then $A_j \lhd_O F_1 \f G_1$, $G'_k \lhd_O F_{k+1} \f
G_{k+1}$ $1 \leq k \leq r-1$, $G_r \approx O(a)$, $G_k \approx G'_k$ $1 \leq k \leq r$, and
$\G\v_{C2_O} t_k:F_k$ $1 \leq k \leq r$. Since $A_j$ is not ending with $O$, then $G_k$ $1 \leq
k \leq r$ is not ending with $O$,  which is impossible since $Cr \approx O(a)$. $\Box$ 
\end{itemize} 

Let $V$ be the set of variables of $\l C$-calculus. \\
Let $P$ be an infinite set of constants called stack constants \footnote{The notion of
stack constants taken from a manuscript of J-L. Krivine}. \\
We define a set of $\l C$-terms $\L CP$ by :
\begin{itemize}
\item[] - If $x \in V$, then $x \in \L CP$ ;
\item[] - If $t \in \L CP$, and $x \in V$, then $\l xt \in \L CP$ ;
\item[] - If $t \in \L CP$, and $u \in \L CP \bigcup P$, then $(t)u \in \L CP$.
\end{itemize}
In other words, $t \in \L CP$ if and only if the stack constants are in argument positions in $t$.\\

Let $\s$ be a function defined on $V \bigcup P$ such that :
\begin{itemize}
\item[] - If $x \in V$, then $\s (x) \in \L CP$ ;
\item[] - If $p \in P$, then $\s (p)=\sou{t}=t_1,...,t_n$, $n \geq 0$, $t_i \in \L CP \bigcup P$
$1 \leq i \leq n$. 
\end{itemize}
We define $\s(t)$ for all $t \in \L CP$ by :
\begin{itemize}
\item[] - $\s ((u)v)=(\s (u))\s (v)$ if $v \not \in P$ ;
\item[] - $\s (\l xu)=\l x \s (u)$ ;
\item[] - $\s ((t)p)=(t)\sou{t}$ if $\s (p)=\sou{t}$.
\end{itemize}
$\s$ is said to be a $P$-substitution.\\

We consider, on the set $\L CP$, the following rules of reduction :
\begin{itemize}
\item[] 1) $(\l xu)tt_1...t_n \f (u[t/x])t_1...t_n$ for all $u,t \in \L CP$ and $t_1,...,t_n \in \L
CP \bigcup P$ ; 
\item[] 2) $(C)tt_1...t_n \f (t)\l x(x)t_1...t_n$ for all $t \in \L CP$ and $t_1,...,t_n \in \L CP
\bigcup P$, and $x$  being $\l$-variable not appearing in $t_1,...,t_n$. 
\end{itemize}

For any  $t,t' \in \L CP$, we shall write $t \rhd_C t'$, if $t'$
is obtained from $t$ by applying these rules finitely many times.  

\begin{lemma} If $t \rhd_C t'$, then $\s (t) \rhd_C \s (t')$ for all $P$-substitution $\s$. 
\end{lemma}
\bf Proof \rm Easy. $\Box$

\begin{lemma} Let $t\in \L CP$ such that the stack constants of $t$ are among
$p_1,...,p_m$. \\ If $t \p_C t'$, and $\G=\G',p_1:O(a_1),...,p_m:O(a_m) \v_{C2_O}
t:\perp$, then $t' \in \L CP$ and $t \rhd_C t'$.  
\end{lemma}
\bf Proof \rm It is enough to do the proof for one step of reduction. We have two cases
: \begin{itemize}
\item[] - If $t=(\l xu)vv_1...v_m$, then, $t'=(u[v/x])v_1...v_m$, $\G,x:F \v_{C2_O}
u:G$, $F$ and $G$ is not ending with $O$, and $\G \v_{C2_O} v:F$. Therefore $u,v \in \L CP$,
and so $t' \in \L CP$ and  $t \rhd_C t'$. 
\item[] - If $t=(C)vv_1...v_m$, then, $t'=(v)\l x(x)v_1...v_m$, and there is a type $A$
which is not ending with $O$ such that $\G \v_{C2_O} v:\neg \neg A$. Therefore $v \in \L CP$,
and so $t' \in \L CP$ and $t \rhd_C t'$. $\Box$ 
\end{itemize} 

\begin{theo} Let $n$ be an integer, $\th_n$ a classical integer of value $n$, and $x,g$ two
distinct variables. 
\begin{itemize}
\item[] - If $n=0$, then for every stack constant $p$, we have : $(\th_n)xgp \p_C
(x)p$. 
\item[] - If $n \not = 0$, then there is $m \geq 1$, and a mapping $IÊ:Ê\{0,...,m\}\f N$,
such that for all distinct stack constants  $p_0,p_1,...,p_m$, we have :
\begin{itemize}
\item[] $(\th_n)xgp_0 \p_C (g)t_1 p_{r_0}$ ;
\item[] $(t_i)p_i \p_C (g)t_{i+1}p_{r_i}$ $1 \leq i \leq m-1$ ; 
\item[] $(t_m)p_m \p_C (x)p_{r_m}$ 
\end{itemize}
where $I(0)=n$, $I(r_m)=0$, and $I(i+1)=I(r_i)-1$ $0 \leq i \leq m-1$.
\end{itemize}  
\end{theo}
\bf Proof \rm  We denote, in this proof, the term $s^i(0)$ by $i$.\\
If $\v_{C2}\th_n:N[n]$, then  $\v_{C2_O} \th_n: [ O(0) \f \perp ], \q y \{ [ O(y) \f \perp ]
\f [O(sy) \f \perp ] \}, O(n) \f \perp $, then  $\G_1= x:O(0) \f \perp, g:\q y \{ [ O(y) \f
\perp ] \f [O(sy) \f \perp ] \}, p_0:O(n) \v_{C2_O} (\th_n)xgp_0:\perp$, therefore, by Lemma 5.1,
$(\th_n)xgp_0$ is $C$-solvable, and three cases may be seen :
\begin{itemize}
\item[]  - If $(\th_n)xgp_0 \p_C (p_0)t_1...t_r$, then $r=0$, and there is a term
$a$, such that $O(a) \approx \perp$. This is impossible. 
\item[] - If $(\th_n)xgp_0 \p_C (x)t_1...t_r$, then $r=1$, and $\G_1 \v_{C2_O} t_1:O(0)$.
Therefore, by Lemma 5.2, $t_1=p_0$, and so $n=0$. 
\item[] - If $(\th_n)xgp_0 \p_C (g)t_1...t_r$, then $r=2$, $\G_1 \v_{C2_O} t_1:O(a) \f \perp$, $\G_1 \v_{C2_O}
t_2:O(s(a'))$, and $a \approx a'$. By Lemma 5.2, we have $t_2=p_0$, and $s(a') \approx n$, then $a
\approx n-1$. Therefore $(\th_n)xgp_0 \p_C (g)t_1p_0$, and $\G_1 \v_{C2_O} t_1:O(n-1) \f \perp$. Let
$I(0)=n$. 
\end{itemize} 
We prove that : if $\G_i=g:\q y \{ [ O(y) \f \perp ] \f [ O(sy) \f \perp ] \} , x:O(0) \f \perp,
p_0:O(I(0)),...., p_i:O(I(i)) \v_{C2_O} (t_i)p_i:\perp$, then : \\
$(t_i)p_i \p_C (g)t_{i+1}p_{r_i}$, and $\G_i \v_{C2_O} t_{i+1}:O(I(r_i)-1) \f \perp$\\
or \\
$(t_i)p_i \p_C (x)p_{r_i}$, and $I(r_i)=0$. \\
$\G_i \v_{C2_O} (t_i)p_i:\perp$, therefore, by Lemma 5.1, $(t_i)p_i$ est $C$-solvable, and three
cases may be seen :  
\begin{itemize}
\item[]  - If $(t_i)p_i \p_C (p_j)u_1...u_r$ $0 \leq j \leq i$, then $r=0$, and there is a
term $a$, such that $O(a) \approx \perp$. This is impossible.  
\item[]  - If $(t_i)p_i \p_C (x)u_1...u_r$, then $r=1$, and $\G_i \v_{C2_O} u_1:O(0)$. Therefore,
by Lemma 5.2, $u_1=p_{r_i}$, and $I(r_i)=0$. 
\item[] - If $(t_i)p_i \p_C (g)u_1...u_r$, then $r=2$, $\G_i \v_{C2_O} u_1:O(a) \f \perp$,
$\G_i \v_{C2_O} u_2:O(s(a'))$, and $a \approx a'$. By Lemma 5.2, we have $u_2=p_{r_i}$, and $s(a')
\approx I(r_i)$, then $a \approx I(r_i)-1$. Therefore $(t_i)p_i \p_C (g)t_{i+1} p_{r_i}$, and
$\G_i \v_{C2_O} t_{i+1}:O(I(r_i)-1) \f \perp$. Let $I(i+1)=I(r_i)-1$.
\end{itemize} 
This construction always terminates. Indeed, if not, the $\l C$-term
$(((\th_n)\l xx)\l xx)p_0$ is not $C$-solvable. This is impossible, since $p_0:\perp \v_{C2}
(((\th_n)\l xx)\l xx)p_0:\perp$. $\Box$  

\begin{corollary} Let $n$ be an integer, $\th_n$ a classical integer of value $n$, and $x,g$
two distinct variables. 
\begin{itemize}
\item[] - If $n=0$, then, for every stack constant $p$, we have : $(\th_n)xgp \rhd_C
(x)p$. 
\item[] - If $n \not = 0$, then there is $m \geq 1$, and a mapping $IÊ:Ê\{0,...,m\}\f N$,
such that for all distinct stack constants  $p_0,p_1,...,p_m$, we have :
\begin{itemize}
\item[] $(\th_n)xgp_0 \rhd_C (g)t_1 p_{r_0}$ ;
\item[] $(t_i)p_i \rhd_C (g)t_{i+1}p_{r_i}$ $1 \leq i \leq m-1$ ; 
\item[] $(t_m)p_m \rhd_C (x)p_{r_m}$ 
\end{itemize}
where $I(0)=n$, $I(r_m)=0$, and $I(i+1)=I(r_i)-1$ $0 \leq i \leq m-1$.
\end{itemize}   
\end{corollary}
\bf Proof \rm We use Lemma 5.4. $\Box$  

\begin{corollary}  Let $n$ be an integer, and $\th_n$ a classical integer of value $n$.
\begin{itemize}
\item[] - If $n=0$, then, for every $\l C-terms$ $a,F,\sou{u}$, we have : $(\th_n)aF\sou{u} \p_C
(a)\sou{u}$.  
\item[] - If $n \not = 0$, then there is $m \geq 1$, and a mapping $IÊ:Ê\{0,...,m\}\f
N$, such that for all $\l C-terms$  $a,F,\sou{u_0},\sou{u_1},...,\sou{u_m}$, we have :
\begin{itemize}
\item[] $(\th_n)aF\sou{u_0} \p_C (g)t_1 \sou{u_{r_0}}$ ; 
\item[] $(t_i)\sou{u_i} \p_C (g)t_{i+1}\sou{u_{r_i}}$ $1 \leq i \leq m-1$ ; 
\item[] $(t_m)\sou{u_m} \p_C (a)\sou{u_{r_m}}$ 
\end{itemize}
where $I(0)=n$, $I(r_m)=0$, and $I(i+1)=I(r_i)-1$ $0 \leq i \leq m-1$.
\end{itemize}    
\end{corollary} 
\bf Proof \rm We use Lemma 5.3. $\Box$  

\subsection{The integers in $M2$}

According to the results of section 4, we can obtain some results concerning the
integers in the system $M2$.

\begin{theo} Let $n$ be an integer. If $\v_{M2} t :N[s^n (0)]$, then, $t \simeq\sb{\b}
\so{n}$. \end{theo}
\bf Proof \rm We use Theorem 4.1. $\Box$ \\

Let $n$ be an integer. By Theorem 4.2, a classical integer of value $n$ is a closed $\l C$-term
$\th_n$ such that $\v_{M2} \th_n :N^C[s^n(0)]$.

\begin{theo}
Let $n$ be an integer, $\th_n$ a classical integer of value $n$, and $x,g$ two 
distinct variables.
\begin{itemize}
\item[] - If $n=0$, then, for every stack constant $p$, we have : $(\th_n)xgp \rhd_C (x)p$.  
\item[] - If $n \not = 0$, then there is $m \geq 1$, and a mapping $IÊ:Ê\{0,...,m\}\f N$, such that for
all distinct stack constants  $p_0,p_1,...,p_m$, we have :
\begin{itemize}
\item[] $(\th_n)xgp_0 \rhd_C (g)t_1 p_{r_0}$ ; 
\item[] $(t_i)p_i \rhd_C (g)t_{i+1}p_{r_i}$ $1 \leq i \leq m-1$ ;  
\item[] $(t_m)p_m \rhd_C (x)p_{r_m}$ 
\end{itemize}
where $I(0)=n$, $I(r_m)=0$, and $I(i+1)=I(r_i)-1$ $0 \leq i \leq m-1$. 
\end{itemize}  
\end{theo}
\bf Proof \rm We use Theorem 4.2. $\Box$

\section{Storage operators}

\subsection{Storage operators for Church integers}

Let $T$ be a closed $\l$-term. We say that $T$ is a storage operator for Church integers if and only if for every
$n \geq 0$, there is a $\l$-term $\t_n \simeq\sb{\b} \so{n}$, such that for every 
$\l$-term $\th_n \simeq\sb{\b} \so{n}$, there is a substitution $\s$, such that $(T)\th_n f \p
(f)\s(\t_n)$.\\

\bf Examples\/ \rm 
If we take :\\ $T_1 = \l n((n)\d)G$ where
$G = \l x\l y(x)\l z(y)(\so{s})z$ and $\d = \l f(f)\so{0}$ \\
$T_2 = \l n\l f(((n)f)F)\so{0}$
where $F = \l x\l y(x)(\so{s})y$, \\
then it is easy to check that : for every $\th_n \simeq\sb{\b} \so{n}$, $(T_i)\th_n f \p
(f)(\so{s})^n \so{0}$ ($i=1$ or $2$) (see [3] and [8]). \\ 
Therefore $T_1$ and $T_2$ are storage operators for Church integers. $\Box$ \\

It is a remarkable fact that we can give simple types to storage operators for Church integers. We
first define the simple G\H{o}del translation $F$* of a formula $F$ : it is obtained by replacing
in the formula $F$, each atomic formula $A$ by $\neg A$. For example :
\begin{center}
$N$*$[x]=\q X \{\neg X(0),\q y(\neg X(y) \f \neg X(sy)) \f \neg X(x) \}$
\end{center}
It is well known that, if $F$ is provable in classical logic, then $F$* is provable in
intuitionistic logic.\\ 

We can check that $\v_{AF2} T_1,T_2 : \q x \{N$*$[x] \f\neg\neg N[x] \}$. And, in general, we have
the following Theorem :

\begin{theo} (see [3] and [10])
If $\v_{AF2} T: \q x\{N$*$[x] \f\neg\neg N[x]\}$, then $T$ is a storage operator for Church integers.
\end{theo}

\subsection{Storage operators for classical integers}

The storage operators play an important role in classical type systems. Indeed, they
can be used to find the value of a classical integer.

\begin{theo} (see [6] and [7])
If $\v_{AF2} T: \q x\{N$*$[x] \f \neg\neg N[x]\}$, then for every $n \geq 0$, there is a
$\l$-term $\t_n \simeq\sb{\b} \so{n}$, such that for every classical integer $\th_n$ of value
$n$, there is a substitution $\s$, such that $(T)\th_n f \p_C (f)\s(\t_n)$.
\end{theo}

\begin{corollary}
If $\v_{AF2} T: \q x\{N$*$[x] \f \neg\neg N[x]\}$, then for every $n \geq 0$ and for every
classical integer $\th_n$ of value $n$, there is a $\l$-term $\t_n$, such that $(T)\th_n \l xx
\p_C \t_n \f\sb{\b} \so{n}$.
\end{corollary}
\bf Proof \rm We use Theorem 6.2. $\Box$ \\

\bf Remark. \rm Theorem 6.2 cannot be generalized for the system $C2$. Indeed, let
$T=\l \n \l f (f) (C)(T_i)\n$ ($i=1$ or $2$). \\

$\n:N$*$[x] , f:\neg N[x] \v_{C2} (T_i)\n:\neg\neg N[x] \Longrightarrow$ \\
$\n:N$*$[x] , f:\neg N[x] \v_{C2} (C)(T_i)\n: N[x] \Longrightarrow$ \\ 
$\n:N$*$[x] , f:\neg N[x] \v_{C2} (f)(C)(T_i)\n: \perp \Longrightarrow$ \\
$\v_{C2} T: \q x\{N$*$[x] \f \neg\neg N[x]\}$ \\

Since for every $\l C$-term $\th$, $(T)\th f \p_C (f) (C)(T_i)\th$, then it is easy to check that
there is not a $\l C$-term $\t_n \simeq\sb{\b} \so{n}$ such that for every classical
integer $\th_n$ of value $n$, there is a substitution $\s$, such that $(T)\th_n f \p_C
(f)\s(\t_n)$. $\Box$ \\ 

We will see that in system $M2$ we have a similar result to Theorem 6.2.\\
 
Let $T$ be a closed $\l C$-term. We say that $T$ is a storage operator for classical integers if
and only if for every $n \geq 0$, there is a $\l C$-term $\t_n \simeq\sb{\b} \so{n}$, such that for
every  classical integers $\th_n$ of value $n$, there is a substitution $\s$, such that $(T)\th_n
f \p_C (f)\s(\t_n)$.

\begin{theo}
If $\v_{M2} T: \q x \{ N^C[x] \f \neg\neg N[x] \}$, then $T$ is a storage operator for classical
integers. 
\end{theo}

The type system $M$ is the subsystem of $M2$ where we only have
propositional variables and constants (predicate variables or predicate symbols of arity 0). So, first
order variable, function symbols, and finite sets of equations are useless. The rules for typed
are $0'$) 1), 2), 3), 6), $6'$), 7) and $7'$) restricted to propositional variables. With each
predicate variable (resp. predicate symbol) $X$, we associate a predicate variable (resp. a
predicate symbol) $X^{\di}$ of $M$ type system. For each formula $A$ of $M2$, we define the
formula $A^{\di}$ of $F_C$ obtained by forgetting in $A$ the first order part. If
$\G=x_1:A_1,...,x_n:A_n$ is a context of $M2$, then we denote by $\G^{\di}$ the context
$x_1:A_1^{\di},...,x_n:A_n^{\di}$ of $M$. We write $\G\v_M t:A$ if $t$ is
typable in $M$ of type $A$ in the context $\G$.\\ 
We have obviously the following property : if $\G \v_{M2} t:A$, then $\G^{\di}
\v_M t:A^{\di}$.\\

Theorem 6.3 is a consequence of the following Theorem.

\begin{theo}
If $\v_M T: N^C \f \neg\neg N$, then for every $n \geq 0$, there is an $m \geq 0$ and a $\l
C$-term $\t_m \simeq\sb{\b} \so{m}$, such that for every classical integer $\th_n$ of value $n$,
there is a substitution $\s$, such that $(T)\th_n f \p_C (f)\s(\t_m)$.  
\end{theo}

Indeed, if $\v_{M2} T: \q x \{ N^C[x] \f \neg\neg N[x] \}$, then $\v_M T:
N^C \f \neg\neg N$. Therefore for every $n \geq 0$, there is an $m \geq
0$ and $\t_m \simeq\sb{\b} \so{m}$, such that for every classical integer $\th_n$ of value $n$,
there is a substitution $\s$, such that $(T)\th_n f \p_C (f)\s(\t_m)$. We have $\v_{M2} \so{n} :
N^C[s^n(0)]$, then $ f:\neg N[s^n(0)] \v_{M2} (T) \so{n} f :\perp$, therefore $ f:\neg
N[s^n(0)]\v_{M2} (f)\so{m} :\perp$ and $\v_{M2} \so{m} : N[s^n(0)]$. Therefore $n=m$.
and $T$ is a storage operator for classical integers. $\Box$\\

In order to prove Theorem 6.4, we shall need some Lemmas.

\begin{lemma} If $\G,\n:N^C\v_M(\n)\sou{d}:\perp$, then $\sou{d}=a,b,d_1,...,d_r$
and there is a classical type $F$, such that : $\G,\n:N^C\v_M a:F$ ;
$\G,\n:N^C\v_M b:F \f F$ ; $F \lhd E_1 \f F_1$, $F_i \lhd E_{i+1} \f F_{i+1}$
$1 \leq i \leq r-1$ ; $F_r \lhd \perp$ ; and $\G,\n:N^C\v_M c_i : E_i$
$1 \leq i \leq r$.  
\end{lemma}
\bf Proof \rm We use Theorem 2.2. $\Box$

\begin{lemma} If $F$ is a classical type and $\G,x:F \v_M (x)\sou{d}:\perp$, then
$\sou{d}=d_1,...,d_r$ ; $F \lhd E_1 \f F_1$ ; $F_i \lhd E_{i+1} \f F_{i+1}$ $1 \leq i \leq r-1$ ; $F_r
\lhd \perp$ ; and $\G,x:F \v_M c_i : E_i$ $1 \leq i \leq r$.  
\end{lemma} 
\bf Proof \rm We use Theorem 2.2. $\Box$

\begin{lemma} Let $t$ be a $\b$-normal $\l C$-term, and $A_1,...,A_n$ a sequence of
classical types.\\  
If $x_1:A_1,...,x_n:A_n \v_M t:N$, then there is an $m \geq 0$ such that $t =
\so{m}$.  
\end{lemma}
\bf Proof \rm  We use Theorems 4.1 and 5.2. $\Box$\\

Let $\n$ and $f$ be two fixed variables.\\ 
We denote by $x_{n,a,b,\sou{c}}$ (where $n$ is an integer, $a,b$ two $\l$-terms, and $\sou{c}$ a
finite sequence of $\l$-terms) a variable which does not appear in $a,b,\sou{c}$.

\begin{theo} 
Let $n$ be an integer. There is an integer $m$ and a finite sequence of head reductions $\{ U_i
\p_C V_i \}_{1\leq i\leq r}$ such that :\\  
1) $U_1 = (T)\n f$ and $V_r = (f)\t_m$ where $\t_m \simeq\sb{\b}\so{m}$ ;\\
2) $V_i = (\n) a b \sou{c}$ or $V_i = (x_{l,a,b,\sou{c}}) \sou{d}$  $0 \leq l \leq n-1$;\\
3) If $V_i = (\n)a b \sou{c}$, then $U_{i+1} = (a)\sou{c}$ if $n=0$ and
$U_{i+1} = ((b)x_{n-1,a,b,\sou{c}})\sou{c}$ if $n \neq 0$ ;\\ 
4) If $V_i = (x_{l,a,b,\sou{c}})\sou{d}$ $0 \leq l \leq n-1$, then $U_{i+1} = (a)\sou{d}$ if $l=0$ and
$U_{i+1} = ((b)x_{l-1,a,b,\sou{d}})\sou{d}$ if $l \neq 0$.
\end{theo}

\bf Proof \rm 
A good context $\G$ is a context of the form $\n:N^C, f:\neg N, x_{n_1,a_1,b_1,\sou{c_1}} : F_1
,..., x_{n_p,a_p,b_p,\sou{c_p}} : F_p$ where $F_i$ is a classical type, $0 \leq n_i \leq n-1$, and $1
\leq i \leq p$ .  \\
We will prove that there is an integer $m$ and a finite sequence of head reductions $\{
U_i \p_C V_i \}_{1\leq i\leq r}$ such that we have 1), 2), 3), 4), and there is a good context $\G$
such that $\G \v_M V_i :\perp$ $1 \leq i \leq r$.\\ 

We have $\v_M T: N^C \f \neg\neg N$, then  $\n:N^C,f:\neg N
\v_M(T)\n f:\perp$, and by Lemmas 6.1 and 6.2, $(T)\n f \p_C V_1$ where $V_1
= (f)\t$ or $V_1 = (\n)ab\sou{c}$.\\  
Assume that we have the head reduction $U_k \p_C V_k$ and $V_k \neq (f)\t$. 
\begin{itemize} 
\item[] 
- If $V_k = (\n)a b \sou{c}$, then, by the induction hypothesis, there is a good context $\G$ such that
$\G \v_M (\n)a b \sou{c} :\perp$. By Lemma 6.1, there is a classical type $F$, such that $\G
\v_Ma:F$ ; $\G \v_Mb:F \f F$ ; $\sou{c}=c_1,...,c_s$ ; $F \lhd E_1 \f F_1$ ;
$F_i \lhd E_{i+1} \f F_{i+1}$ $1 \leq i \leq s-1$ ; $F_s \lhd \perp$ ; and $\G \v_Mc_i :
E_i$ $1 \leq i \leq s$.   
\begin{itemize}
\item[]
- If $n=0$, let  $U_{k+1} = (a)\sou{c}$. We have $\G \v_M U_{k+1}:\perp$. 
\item[]
- If $n \neq 0$, let $U_{k+1} = ((b)x_{n-1,a,b,\sou{c}})\sou{c}$. The variable
$x_{n-1,a,b,\sou{c}}$ is not used before. Indeed, if it is, we check easily that the $\l C$-term
$(T)\so{n} f$ is not solvable; but that is impossible because $f:\neg N \v_M(T)\so{n} f :\perp$.
Therefore $\G' = \G ,x_{n-1,a,b,\sou{c}}:F$ is a good context and $\G' \v_M U_{k+1} :\perp$. 
\end{itemize}  
\end{itemize}
\begin{itemize}
\item[] - If $V_k = (x_{l,a,b,\sou{c}}) \sou{d}$, then, by the induction hypothesis, there is a
good context $\G$ such that $\G \v_M(x_{l,a,b,\sou{c}}) \sou{d} :\perp$. 
Then there is a classical type $F$ such that $x_{l,a,b,\sou{c}} : F$ is in the context $\G$. By
Lemma 6.2, $\sou{c}=d_1,...,d_s$ ; $F \lhd E_1 \f F_1$ ; $F_i \lhd E_{i+1} \f F_{i+1}$ $1 \leq i
\leq s-1$ ; $F_s \lhd \perp$ ; and $\G \v_Mc_i : E_i$ $1 \leq i \leq s$.  
\begin{itemize}
\item[] 
- If $l=0$, let $U_{k+1}=(a)\sou{c}$. We have $\G \v_M U_{k+1}:\perp$.
\item[] 
- If $l \neq 0$,  Let $U_{k+1} = ((b)x_{l-1,a,b,\sou{d}})\sou{d}$. The variable
$x_{l-1,a,b,\sou{d}}$ is not used before. Indeed, if it is, we check that the $\l C$-term
$(T)\so{n} f$ is not solvable; but this is impossible because $f:\neg N \v_M (T)\so{n} f :\perp$.
Then $\G' = \G ,x_{l-1,a,b,\sou{d}}:F$ is a good context and $\G' \v_M U_{k+1} :\perp$. 
\end{itemize}
\end{itemize}
Therefore there is a good context $\G'$ such that $\G' \v_M U_{k+1} :\perp$. Then, by
Lemmas 6.1 and 6.2, $U_{k+1} \p_C V_{k+1}$ where  $V_{k+1} = (f)\t$ or $V_{k+1} =
(\n)ab\sou{c}$ or  $V_{k+1} = (x_{l,a,b,\sou{c}})\sou{d}$ $0 \leq l \leq n-1$. \\ 
This construction always terminates. Indeed, if not, we check that the $\l C$-term $(T)\so{n} f$ is
not solvable; but this is impossible because $f:\neg N \v_M (T)\so{n} f :\perp$. \\ 
Therefore there is $r \geq 0$ and a good context $\G$ such that $\G \v_M V_{r} = (f)\t
:\perp$, and $\G \v_M \t :N$. Therefore, by Lemma 6.3, there is an $m \geq 0$ such that $\t
\simeq\sb{\b}\so{m}$. $\Box$ \\

Let $T$ be a $\l C$-term such that $\v_M T: N^C \f \neg\neg N$. 
By Theorem 6.5, there is an integer $s$ and a finite sequence of head reductions $\{ U_i \p_C V_i
\}_{1\leq i\leq r}$ such that :\\  
1) $U_1 = (T)\n f$ and $V_r = (f)\t_s$ where $\t_s \simeq\sb{\b}\so{s}$;\\
2) $V_i = (\n) a b \sou{c}$ or $V_i = (x_{l,a,b,\sou{c}}) \sou{d}$  $0 \leq l \leq n-1$;\\
3) If $V_i = (\n)a b \sou{c}$, then $U_{i+1} = (a)\sou{c}$ if $n=0$ and
$U_{i+1} = ((b)x_{n-1,a,b,\sou{c}})\sou{c}$ if $n \neq 0$ ;\\ 
4) If $V_i = (x_{l,a,b,\sou{c}})\sou{d}$ $0 \leq l \leq n-1$, then $U_{i+1} = (a)\sou{d}$ if $l=0$ and
$U_{i+1} = ((b)x_{l-1,a,b,\sou{d}})\sou{d}$ if $l \neq 0$.\\

Let $\th_n$ be a classical integer of value $n$, and $x,g$ two distinct variables. By Theorem
5.6 we have :\\ 
If $n=0$, then for every stack constant $p$, we have : $(\th_n)xgp \rhd_C (x)p$. \\
If $n \not = 0$, then there is $m \geq 1$, and a mapping $IÊ:Ê\{0,...,m\}\f N$, such that for
all distinct stack constants  $p_0,p_1,...,p_m$, we have :\\
$(\th_n)xgp_0 \rhd_C (g)t_1 p_{r_0}$ ;\\
$(t_i)p_i \rhd_C (g)t_{i+1}p_{r_i}$ $1 \leq i \leq m-1$ ; \\
$(t_m)p_m \rhd_C (x)p_{r_m}$ \\
where $I(0)=n$, $I(r_m)=0$, and $I(i+1)=I(r_i)-1$ $0 \leq i \leq m-1$.\\

\begin{lemma}If $n = 0$, then $(T)\th_n f \p_C (f)\t[\th_n / \n]$.
\end{lemma}
\bf Proof \rm We prove by induction that for every $1 \leq i \leq r$, we have $(T)\th_n f \p_C V_i
[\th_n / \n]$.\\ For $i=1$, $(T)\th_n f = \{ (T)\n f \} [\th_n / \n] = U_1 [\th_n / \n] \p_C V_1
[\th_n / \n]$.\\   Assume it is true for $i$, and prove it for $i+1$.\\
$(T)\th_n f \p_C V_i [\th_n / \n] = \{ (\n)ab \sou{c} \} [\th_n / \n] = \{ (\th_n)ab\sou{c} \}
[\th_n / \n] = \{ (\th_n)xgp \} [a / x , b / g , \sou{c} / p ] [\th_n / \n]$. Since $(\th_n)xgp
\p_C (x)p$, then $(T)\th_n f \p_C \{ (a)\sou{c} \} [\th_n / \n] = U_{i+1} [\th_n / \n] \p_C
V_{i+1}[\th_n / \n]$.\\
So, for $i=r$, we have $(T)\th_n f \p_C V_r [\th_n / \n] = \{ (f)\t \} [\th_n / \n] = (f)\t[\th_n /
\n]$. $\Box$ \\

We assume now that $n \geq 1$. \\

A $k-\l C$-term is a $\l C$-term of the forme
$V_k[\t_1 / y_1 ] ... [\t_p / y_p ] [\th_n / \n ]$ such that : \\
- $Fv(V_k) \subseteq \{ \n , f , y_1 , ... , y_p \}$ \\
- for every $1 \leq i \leq p$, $y_i = x_{n_i,a_i,b_i,\sou{c}_i}$ and $\t_i =t_{m_i}
[a_i / x , b_i / g , \sou{d_0} / p_0 ,...,\sou{d_{m_i-1}} / p_{m_i-1}]$ where
$I(m_i)=n_i$ \\
- for every $0 \leq k \leq m_i-1$, there is $1 \leq l \leq r$ such that $U_l =
(a_i)\sou{d_k}$ if $I(k)=0$ and $U_r = (b_i)x_{I(k)-1,a_i,b_i,\sou{d_k}}\sou{d_k}$ if
$I(k) > 0$.\\

To simplify, a $k- \l C$-term is denoted by $V_k []$.

\begin{lemma} Let $1 \leq i \leq r-1$ and $V_i []$ an $i-\l C$-term. If $(T)\th_n f \p_C V_i []$,
then there is $1 \leq j \leq r$ and a $j-\l C$-term $V_j []$ such that $V_j [] \p_C V_j []$ and
either $V_i [] \not = V_j []$ or $i < j$ 
\end{lemma} 
\bf Proof \rm There are only two possibilities. 1) $V_i =
(\n)ab\sou{c}$ ;  2) $V_i = (x_{\a,a,b,\sou{c}})\sou{d}$. \\
We now examine each of this cases.\\
1) If $V_i = (\n)ab\sou{c}$, then $V_i [] = \{ (\th_n)ab\sou{c} \}[] = \{ (\th_n)xgp_0 \} [a / x
, b / g ,\sou{c} / p_0] []$. Since $(\th_n)xgp_0 \rhd_C (g)t_1 p_{r_0}=(g)t_1 p_0$, then $V_i []
\p_C \{ (b)t_1[a / x, b / g , \sou{c} / p_0 ] \sou{c} \}  [] =$ \\
$\{ (b) x_{n-1,a,b,\sou{c}} \sou{c}  \} [t_1[a / x, b / g , \sou{c} / p_0 / x_{n-1,a,b,\sou{c}}]
[] = U_{i+1} [] \p_c V_{i+1} []$. Let $j=i+1$. We have $i < j$ and $I(1)=I(r_0)-1=I(0)-1=n-1$.\\
2) If $V_i = (x_{\a,a,b,\sou{c}})\sou{d}$, then $V_i []=\{ (t_{\b}[a / x, b / g ,  \sou{d_0} / p_0
,...,\sou{d_{\b-1}} / p_{\b-1}])\sou{d} \} []$ where $I(\b)=\a$.\\
If $I(\b)=\a \not =0$, then
$U_{i+1}=(b)x_{\a-1,a,b,\sou{d}}\sou{d}=(b)x_{I(\b)-1,a,b,\sou{d}}\sou{d}$, and if
$I(\b)=\a \not =0$, then $U_{i+1}= (a)\sou{d}$.\\
We consider the following two cases.
\begin{itemize}
\item[] - If $\b \leq m$, then $(t_{\b})p_{\b} \rhd_C (g)t_{{\b}+1}p_{r_{\b}}$, so that \\
$V_i [] \p_C \{ (g)t_{{\b}+1}p_{r_{\b}} \}[a / x, b / g ,  \sou{d_0} / p_0
,...,\sou{d_{\b-1}} / p_{\b-1} , \sou{d} / p_{\b}] []$ = \\
$\{ (b)t_{{\b}+1} \sou{d_{r_{\b}}} \}[a / x, b / g , \sou{d_0} / p_0 ,...,\sou{d_{\b-1}} /
p_{\b-1},\sou{d} / p_{\b}] []$. \\ 
Since $\b \not = m$, then $I(r_\b) \not = 0$.
By the hypothesis there is $1 \leq j \leq r$ such that $U_j =
(b)x_{I(r_{\b})-1,a,b,\sou{d_{r_{\b}}}}\sou{d_{r_{\b}}}$. Therefore \\
$V_i [] \p_C U_j [t_{{\b}+1}[a / x, b / g ,  \sou{d_0} / p_0 ,...,\sou{d_{\b-1}} /
p_{\b-1},\sou{d} / p_{\b}]/ x_{I(r_{\b})-1,a,b,\sou{d_{r_{\b}}}}][] =
U_j [] \p_C V_j []$.\\  
If $V_i [] = V_j []$, then the head $C$-reduction $(t_{\b})p_{\b} \rhd_C
(g)t_{{\b}+1}p_{r_{\b}}$ must be an identity, in other words $(t_{\b})p_{\b} =
(g)t_{{\b}+1}p_{r_{\b}}$ and therefore $\b = r_{\b}$. And so $j = i+1 > i$. 
\item[] - If $\b = m$, then $(t_{\b})p_{\b} = (t_m)p_m \rhd_C (x)p_{r_m}$, so that \\
$V_i [] \p_C \{ (x)p_{r_m} \}[a / x, b / g ,  \sou{d_0} / p_0 ,...,\sou{d_{m-1}} / p_{m-1}][] =
(a)t_{r_m} \}[]$.\\ Since $I(r_m)=0$, then by the hypothesis there is $1 \leq j \leq r$
such that $U_j = (a)t_{r_m}$. Therefore  $V_i [] \p_C U_j [] \p_C V_j []$.\\
If $V_i [] = V_j []$, then the head $C$-reduction $(t_m)p_m \rhd_C (x)p_{r_m}$
must be an identity, in other words $(t_m)p_m \rhd_C (x)p_{r_m}$ and therefore $m
= r_m$. And so $j = i+1 > i$. $\Box$
\end{itemize}

\begin{corollary}
There is a substitution $\s$ such that $(T)\th_n f \p_C (f)\s(\t)$.
\end{corollary}
\bf Proof \rm $(T)\th_n f = \{ (T)\n f \} [\th_n / \n] = U_1 [\th_n / \n]
\p_C V_1 [\th_n / \n]$. By Lemma 6.5 we obtaine a sequence $V_{i_1} []$ , $V_{i_2} []$ , ... , 
$V_{i_k} []$ , ... such that $(T)\th_n f \p_C V_{i_s} []$ and if $V_{i_s} [] \not = V_{i_{s+1}} []$
then $i_s \leq i_{s+1}$. This sequence is necessarily finite, indeed $f:\neg N \v_M (T)\th_n f
:\perp$. If $V_{i_s} [] = V_{i_{s+1}} [] = ... =  V_{i_{s+\a}} []$, then  $i_s
< i_{s+1} < ... < i_{s+\a}$ and $\a \leq r$. Therefore there is $s$ such
that $V_{i_s} = (f)\t$, then $(T)\th_n f \p_C  V_{i_s} [] = \{ (f)\t \} [] = (f)\t[]$. $\Box$ \\

Then, by Lemma 6.4 and Corollary 6.2, $T$ is a storage operator for classical integers.

\subsection{General Theorem}

In this subsection, we give (without proof) a generalization of Theorem 6.3.\\

Let $T$ be a closed $\l C$-term, and $D,E$ two closed types of $AF2$ type system. We say that $T$ is
a storage operator for the pair of types $(D,E)$ iff for every $\l$-term $\v_{AF2} t:D$, there is
$\l$-term $\t'_t$ and $\l C$-term $\t_t$, such that $\t'_t \simeq\sb{\b} \t_t$, $\v_{AF2}
\t'_t:E$, and for every $\v_{C2} \th_t:D$, there is a substitution $\s$, such that $(T)\th_t f
\p_C (f)\s(\t_t)$.

\begin{theo} Let $D,E$ two $\q$-positive closed types of $AF2$ type system, such that $E$ does not
contain $\perp$. If $\v_{M2} T: D^C \f \neg\neg E$, then $T$ is a storage operator for the
pair $(D,E)$.   
\end{theo}

\section{Operational characterization of $\l C$-terms of type $\q X_C \{\perp \f X_C \}$ and
$\q X_C \{ \neg \neg X_C \f X_C \}$}
 
Let \bf A \rm (for Abort) the $\l C$-term $\l x(C)\l yx$.\\ 

\bf Behaviour of \bf A \rm : \rm 
\begin{center}
(\bf A \rm)$tt_1...t_n \p_C ((C)\l yt) t_1...t_n \p_C (\l yt)\l x(x)t_1...t_n \p_C t$.
\end{center}
\bf Typing of \bf A : \rm
\begin{center}
$x:\perp \v_{M2} \l yx:\neg \neg X_C \Longrightarrow x:\perp \v_{M2} (C)\l yx:X_C
\Longrightarrow \v_{M2}$ \bf A \rm $:\q X_C \{ \perp \f X_C \}$ 
\end{center}

\begin{theo} If $\v_{M2} T:\q X_C \{ \perp \f X_C \}$, then for every integer $n$, and
for all $\l C-terms$ $t,t_1,...,t_n$, $(T)t t_1...t_n \p_C t$. 
\end{theo}
\bf Proof. \rm Let $O_1,...,O_n$ be new predicate symbols of arity 0 different from
$\perp$. Let $A=O_1,...,O_n \f \perp$. If $\v_{M2} T:\q X_C \{ \perp \f X_C \}$, then
$\v_{M2} T:\perp \f A$, and $\G = x:\perp,x_1:O_1,...,x_n:O_n \v_{M2} (T)xx_1...x_n:\perp$.
Therefore $(T)xx_1...x_n \p_C (f)u_1...u_r$ and $\G \v_{M2} (f)u_1...u_r:\perp$.
\begin{itemize} 
\item[] - If $f=x_i$ $1 \leq i \leq n$, then $r=0$, and $O_i=\perp$. A contradiction. 
\item[] - If $f=x$, then $r=0$, and $(T)xx_1...x_n \p_C x$, therefore, for every integer $n$, and
for all $\l C$-terms $t,t_1,...,t_n$, $(T)t t_1...t_n \p_C t$. $\Box$ 
\end{itemize} 

The constant $C$ satisfies the following relations :\\
$(C)t t_1...t_n \p_C  (t)U$ and\\
$(U)y \p_C (y)t_1...t_n$ where $y$ is a new variable. \\ 

Let $C'=\l x(C)\l d(x)\l y(x)\l z(d)y$. \\

$x:\neg \neg X_C,y:X_C,z:X_C,d:\neg X_C \v_{M2} (d)y:\perp \Longrightarrow$ \\
$x:\neg \neg X_C,y:X_C,d:\neg X_C \v_{M2} (x)\l z(d)y:\perp \Longrightarrow$ \\ 
$x:\neg \neg X_C,d:\neg X_C \v_{M2} (x)\l y(x)\l z(d)y:\perp \Longrightarrow$ \\
$x:\neg \neg X_C \v_{M2} (C)\l d(x)\l y(x)\l z(d)y:X_C \Longrightarrow$ \\
$\v_{M2} C': \q X_C \{\neg \neg X_C \f X_C \}$. \\ 

The $\l C$-term $C'$ satisfies the following
relations :\\ 
$(C')t t_1...t_n \p_C (t)U$,\\ 
$(U)y \p_C (t)V$, and \\
$(V)z \p_C (y)t_1...t_n$ where $y,z$ are new variables. \\

In general, we have the following characterization.

\begin{theo}  If $\v_{M2} T: \q X_C \{\neg \neg X_C \f X_C \}$, then there is an
integer $m$, such that, for every integer $n$, and for all $\l C$-terms
$t,t_1,...,t_n$ :
\begin{itemize} 
\item[] $(T)t t_1...t_n \p_C (t)V_1$, 
\item[] $(V_i)y_i \p_C (t)V_{i+1}$ $1 \leq i\leq m-1$, and 
\item[] $(V_m)y_m \p_C (y_i)t_1...t_n$ where $y_1,...,y_m$ are new variables.
\end{itemize}  
\end{theo} 

\bf Proof \rm Let $O$ be a new predicate symbol of arity 0 different from $\perp$.
We define as in section 3, the system $M2_O$. And we check easily that this system has the
same results as Lemmas 5.1, 5.2, 5.3 and 5.4.\\ 
Let $p$ be a stack constant and $A=O \f \perp$. If $\v_{M2} T: \q X_C \{\neg \neg X_C \f X_C \}$, 
then $\v_{M2_O} T: \neg \neg A \f A$, and $\G=x:\neg \neg A, p:O
\v_{M2_O} (T)xp:\perp$. Therefore $(T)xp \p_C (f)u_1...u_r$, and $\G \v_{M2_O}
(f)u_1...u_r:\perp$.
\begin{itemize} 
\item[] - If $f=p$, then $r=0$, and $O=\perp$. A contradiction. 
\item[] - If $f=x$, then, $(T)xp \rhd_C (x)U_1$, and $\G \v_{M2_O} U_1:\neg A$.
\end{itemize}  
We prove (by induction) that if $\G,y_1:A,...,y_{i-1}:A \v_{M2_O} U_i:\neg A$, then
[$(U_i)y_i \rhd_C (x)U_{i+1}$, and $\G,y_1:A,...,y_i:A \v_{M2_O} U_{i+1}:\neg
A$] or [$(U_i)y_i \rhd_C (y_j)p$ $1 \leq j \leq i$]. \\
The sequence $(U_i)_{i \geq 0}$ is not infinite. Indeed, if it is, the $\l C$-term $((T)\l x(x)z)p$
is not $C$-solvable; but this is impossible, because $z:A,p:O \v_{M2} ((T)\l x(x)z)p:\perp$.
\\ To obtain the Theorem, we replace the constant $p$ by the sequence $\sou{t}=t_1,...,t_n$ and we
put $V_i = U_i [\sou{t} / p]$. $\Box$

\section{The $\l\m$-calculus}

In this section, we give a similar version to Theorem 6.3 in the M. Parigot's $\l \m$-calculus.

\subsection{Pure and typed $\l\m$-calculus}

$\l\m$-calculus has two distinct alphabets of variables : the set of $\l$-variables
$x,y,z,...$, and the set of $\m$-variables $\a,\b,\g$,.... Terms are defined
by the following grammar :
\begin{center}
$t$ $:=$ $x$ $\mid$ $\l xt$ $\mid$ $(t)t$ $\mid$  $\m\a[\b]t$
\end{center}

Terms of $\l\m$-calculus are called $\l\m$-terms.\\  

The reduction relation of $\l\m$-calculus is induced by fives different notions of
reduction :\\ 

\bf The computation rules \/ \rm

\begin{itemize} 
\item[] ($C_1$) $(\l xu)v \f u[v/x]$ 
\item[] ($C_2$) $(\m\a u)v \f \m\a u[v/$*$\a]$  
\item[] where $u[\sou{v}/$*$\a]$ is obtained from $u$ by replacing inductively each
subterm of the form $[\a]w$ by $[\a](w)\sou{v}$.  
\end{itemize} 

\bf The simplification rules \/ \rm

\begin{itemize} 
\item[] ($S_1$)  $[\a]\m\b u \f  u[\a/\b]$  
\item[] ($S_2$)  $\m\a [\a]u \f u$, if $\a$ has no free occurence in $u$   
\item[] ($S_3$)  $\m\a u \f \l x \m\a u[x/$*$\a]$, if $u$ contains a subterm of the
form $[\a]\l yw$.  
\end{itemize}

\begin{theo} (see [18]) In $\l\m$-calculus, reduction is confluent.
\end{theo} 

The notation $u \p_{\m} v$ means that $v$ is obtained from $u$ by some head reductions.\\ 
The head equivalence relation is denoted by : $u \sim_{\m} v$ if and only if there is a $w$, 
such that $u \p_{\m} w$ and $v \p_{\m} w$.\\

Proofs are written in a natural deduction system with several conclusions, presented
with sequents. One deals with sequents such that :\\  
- Formulas to the left of $\v$ are labelled with $\l$-variables ;\\  
- Formulas to the right of $\v$ are labelled with $\m$-variables, except one formula
which is labelled with a $\l\m$-term ;\\  
- Distinct formulas never have the same label.\\
  
The right and the left parts of the sequents are considered as sets and therefore
contraction of formulas is done implicitly. \\
  
Let $t$ be a $\l\m$-term, $A$ a type, $\G = x_1:A_1,...,x_n:A_n$,
and $\D = \a_1:B_1,...,\a_m:B_m$. We define by means of the following
rules the notion "$t$ is of type $A$ in $\G$ and $\D$". This notion is denoted by
$\G\v_{FD2} t:A,\D$. 
\begin{itemize} 
\item[] The rules (1),...,(8) of $AF2$ type system.
\item[] (9) If $\G\v_{FD2} t:A,\b:B,\D$, then $\G\v_{FD2}\m\b [\a]t:B,\a:A,\D$.    
\end{itemize} 

Weakenings are included in the rules (2) and (9).\\

As in typed $\l$-calculus on can define $\neg A$ as $  \f \perp$ and use the previous rules with
the following special interpretation of naming for $\perp$ : for $\a$ a $\m$-variable, $\a :
\perp$ is not mentioned.\\

\bf Example \rm Let \bf C \rm =$\l x \m \a [\ph](x)\l y \m \b [\a] y$.\\
$x:\neg \neg X,y:X \v_{FD2} y:X \Longrightarrow$ \\
$x:\neg \neg X,y:X \v_{FD2} \m \b [\a]y:\perp,\a:X \Longrightarrow$ \\ 
$x:\neg \neg X \v_{FD2} \l y \m \b [\a]y:\neg X,\a:X \Longrightarrow$ \\ 
$x:\neg \neg X \v_{FD2} \m \a [\ph] (x) \l y \b [\a]y: X \Longrightarrow$ \\ 
$\v_{FD2}$\bf C \rm $: \q X \{\neg \neg X \f X \}$.  

\begin{theo} (see [18] and [20]) The $FD2$ type system has the following properties :\\
1) Type is preserved during reduction.\\
2) Typable $\l\m$-terms are strongly normalizable. 
\end{theo}

\subsection{Classical integers}

Let $n$ be an integer. A classical integer of value $n$ is a closed $\l\m$-term $\th_n$ such that
$\v_{FD2} \th_n :N[s^n(0)]$.  \\

Let $x$ and $f$ fixed variables, and $N_{x,f}$ be the set of $\l\m$-terms defined by the
following grammar :
\begin{center}
$u$ $:=$ $x$ $\mid$ $(f)u$ $\mid$ $\m\a[\b]x$ $\mid$ $\m\a[\b]u$
\end{center}

We define, for each $u \in N\sb{x,f}$ the set $rep(u)$, which is intuitively the set
of integers potentially repesented by $u$ :
\begin{itemize}
\item[] - $rep(x) = \{ 0 \}$ 
\item[] - $rep((f)u) = \{ n+1$ if $n \in rep(u) \}$ 
\item[] - $rep(\m\a[\b]u)= \bigcap rep(v)$ for each subterm $[\a]v$ of $[\b]u$
\end{itemize}

The following Theorem characterizes the normal forms of classical integers.

\begin{theo} (see [19]) The normal classical integers of value $n$ are exactly the
$\l\m$-terms of the form $\l$x$\l$fu with u$\in$ $N_{x,f}$ without free $\m$-variable
and such that rep(u)=$\{n\}$.
\end{theo}

\subsection{General Theorem}

In order to define, in this framework, the equivalent of system $M2$, the
demonstration of $\neg \neg A \f A$ should not be allowed for all formulas A, and
thus we should prevent the occurrence of some formulas on the right. Thus the following
definition. \\ 

Let $t$ be a $\l\m$-term, $A$ a type, $\G = x_1:A_1,...,x_n:A_n$,
and $\D = \a_1:B_1,...,\a_m:B_m$ where $B_i$ $1 \leq i \leq m$ is a classical
type. We define by means of the following rules the notion "$t$ is of type $A$ in
$\G$ and $\D$", this notion is denoted by $\G\v_{M2} t:A,\D$.

\begin{itemize}
\item[] The rules of $DL2$ type system. 
\item[] (6$'$) If $\G\v t:A, \D$, and $X_C$ has no free occurence in $\G$, then $\G\v t: \q
X_C A, \D$.  
\item[] (7$'$) If $\G\v t: \q X_C A, \D$, and $G$ is a classical type, then
$\G\v t:A[G/ X_C], \D$.  
\end{itemize} 

Let $T$ be a closed $\l\m$-term. We say that $T$ is a storage operator for classical
integers if and only if for every $n \geq 0$, there is $\l \m$-term $\t_n \simeq\sb{\b} \so{n}$,
such that for every  classical integers $\th_n$ of value $n$, there is a substitution
$\s$, such that $(T)\th_n f  \sim_{\m} \m \a [\a] (f)\s(\t_n)$.

\begin{theo}
If $\v_{M2} T: \q x \{ N^C[x] \f \neg\neg N[x] \}$, then $T$ is a storage operator for
classical integers. 
\end{theo}

\end{document}